\documentclass [11pt,reqno]{amsart}
\usepackage {amsmath,amssymb,verbatim,geometry}
\usepackage[colorlinks]{hyperref}
\usepackage[all]{xy}
\newif\ifpdf
\ifpdf
\usepackage[pdftex]{hyperref}
\else
\fi

\numberwithin{equation}{section}       

\newtheorem{prop} {Proposition} [section]
\newtheorem{thm}[prop] {Theorem} 
 
\newtheorem{defi}[prop] {Definition}
\newtheorem{lem}[prop] {Lemma}

\newtheorem{prop-def}[prop]{Proposition-Definition}

\newtheorem{exa}[prop]{Example}
\newtheorem*{thmA}{Theorem A} 
 
\newtheorem*{thmB}{Theorem B}

\theoremstyle{remark}
\newtheorem{rem}[prop]{Remark}

\newcommand{\C}{{\mathbb{C}}}

\newcommand{\R}{{\mathbb{R}}}
\newcommand{\Rn}{{\mathbb R^n}}

\newcommand{\cH}{{\mathcal{H}}}

\newcommand{\psh}{{\mathrm{PSH}}}

\newcommand{\Om}{\operatorname{\Omega}}

\newcommand{\lbr}{{\lbrace}}
\newcommand{\rbr}{{\rbrace}}

\title[Regularity of geodesics]{Regularity of geodesics in the spaces of convex and plurisubharmonic functions III}

\date{\today}

\author{Soufian Abja, S\l awomir Dinew}

\keywords{Monge-Amp\`ere equation, regularity}
\address{École Normale Supérieure, Moulay Ismail University, 3104 Toulal, Meknès, Morocco}
\email{s.abja@umi.ac.ma}
\address{Institute of Mathematics, Jagiellonian University, ul Lojasiewicza 6, 30-348 Krak\'ow, Poland}
\email{Slawomir.Dinew@im.uj.edu.pl }
\subjclass[2010]{Primary: 35J96, secondary: 35J70}

\begin{document}
	\begin{abstract}
		In this note, we study the regularity of solutions to the geodesic equation in the space of plurisubharmonic functions in a general Reinhardt domain. Our findings resolve a problem left open from our previous note \cite{AD26}.
	\end{abstract}
	
	\maketitle

	\section*{Introduction} This is a follow-up note of our previous manuscripts \cite{AD21}, \cite{AD26}. Recall that, modelling a K\"ahler geometric construction Rashkovskii \cite{Ras17}, \cite{Ra17} and independently the first named author \cite{Abj19} studied the geodesic equation in the local space of plurisubharmonic functions.

	Given a bounded hyperconvex domain $\Om\subset\C^n$, $n\geq 1$ let $\psh(\Om)$ stand for the space of plurisubharmonic functions on $\Omega$ and let $\cH\subset\psh(\Om)$ be the subset of smooth strictly plurisubharmonic functions vanishing on the boundary and having  finite $C^2$-norm. Given two functions $\varphi_0,\ \varphi_1\in\cH$.
The geodesic joining $\varphi_0$ and $\varphi_1$ is the  solution (if any) of the following Dirichlet problem (see \cite{AD21})
\begin{equation}\label{MAcomplex}
	\begin{cases}
		\phi\in \psh(\Omega\times A)\cap C(\overline{\Omega\times A});\\
		(dd^c_{z,\zeta}\phi)^{n+1}=0 &{\rm in}\  \hbox{$\Omega\times A$}; \\
		\phi=\varphi_0 &{\rm in}\ \hbox{$\Omega\times \{|\zeta|=1\}$}; \\
		\phi=\varphi_1 &{\rm in}\ \hbox{$\Omega\times \{|\zeta|=e\}$}; \\
		\phi=0 &{\rm in}\ \hbox{$\partial \Omega\times A$,}
		\tag{DMA}
	\end{cases}
\end{equation}
where  $A=\{z\in\C:  1<|z|<e\}$ is the annulus in $\C$ and  $(dd^c_{z,\zeta}\cdot)^n$ is the complex Monge-Amp\`ere operator in joint variables $(z,\zeta)$ which, for smooth input $\phi$, is the determinant of the complex Hessian of $\phi$.

It is a classical fact in pluripotential theory that (\ref{MAcomplex}) admits unique weak solution which has to be rotationally invariant with respect to the annular variable $\zeta$ i.e.  $\phi(z,\zeta)=\varphi_t(z)$ with $t=\log|\zeta|$.

Very little is known about the regularity of the weak solution $\phi$ in general. If $\Omega$ has Lipschitz boundary then $\psi$ is globally Lipschitz as the standard Walsh technique (\cite{Wal68})  shows (see \cite{AD21} for details). Note however that the product domain $\Om\times A$ cannot have better than Lipschitz boundary because of the corners and thus the higher order regularity requires much subtler tools. We refer to the recent paper \cite{DR26}, which even treats more general boundary conditions.

On the bright side much more can be said if $\Omega$ is assumed to have certain symmetries. In \cite{AD21}, \cite{AD26}, we investigated the regularity of $\phi$ in complete Reinhardt domains. In such a setting (\ref{MAcomplex}) admits a reformulation to a real Monge-Amp\`ere problem on the logarithmic image of $\Omega$. The next theorem recaptures the findings from \cite{AD21} and \cite{AD26}:
\begin{thmA}\label{complgeod2}
 Let $\Omega$ be a smoothly bounded strictly pseudoconvex complete Reinhardt domain in $\mathbb C^n$. Suppose that $\phi$ is a weak geodesic solving the problem (\ref{MAcomplex}). If $\varphi_0,\varphi_1\in\mathcal{\tilde{H}}$ are toric in the space variables i.e.
	for all $z,z'\in \Omega$ satisfying $|z|=|z'|$ one has $\varphi_j(z)=\varphi_j(z'),\ j=0,1$ then $\phi$ is $C^{1,1}$ away from the corner $\partial\Omega\times\partial A.$
\end{thmA}

The question that remained open (see \cite{AD26}, Question 1) is whether the completeness assumption on $\Omega$ is needed. This technical assumption was used in \cite{AD26} to control the geometry of the line segments on $U\times[0,1]$ (with $U$ being the logarithmic image of $\Omega$) along which the convex function $u(x,t)=\phi(e^{x_1},\cdots,e^{x_n},e^t)$ is affine. More specifically, given any point $(x,t)\in U\times[0,1]$ there are unique points $\xi=\xi(x,t),\eta=\eta(x,t)\in U$ so that the line segment $l$ in $U\times[0,1]$ with endpoints $(\xi,0)$ and $(\eta,1)$ passes through $(x,t)$ and the convex function $u(x,t):=\phi(e^{x_1},\cdots,e^{x_n},e^t)$ is affine along $l$. The open problem can be described as follows:

{\bf Question}: {\it Let $\Omega$ be smoothly bounded strictly pseudoconvex Reinhardt domain (not necessarily complete).  Is the vector valued $\xi-\eta$ function bounded?}

As explained in \cite{AD26} this technical property is crucial for the control of the $C^{1,1}$ norm of the solution near the coordinate axes.

The goal of this note is to answer affirmatively this question. Surprisingly (for us) the geometric properties of $\Omega$ alone yield the argument regardless of $\varphi_0$ and $\varphi_1$. As a corollary, we obtain the following main theorem of  this note.
\begin{thmB}\label{complgeod3}
 Let $\Omega$ be a smoothly bounded strictly pseudoconvex Reinhardt domain in $\mathbb C^n$. Suppose that $\phi$ is a weak geodesic solving the problem (\ref{MAcomplex}). If $\varphi_0,\varphi_1\in \tilde{\mathcal H}$ are toric in the space variables i.e.
	for all $z,z'\in \Omega$ satisfying $|z|=|z'|$ one has $\varphi_j(z)=\varphi_j(z'),\ j=0,1$ then $\phi$ is $C^{1,1}$ away from the corner $\partial\Omega\times\partial A.$
\end{thmB}

{\bf Acknowledgements}. The first name author was partially supported by the Qatar National Research Fund, project NPRP-BSRA01-0309-210004. The second named author is  supported in part by National Science Center of Poland grant no. 2025/57/B/ST1/00367.

\section{preliminaries}

To begin with we recall the notion of Reinhardt domains:

\begin{defi}\label{Reinhardt}
	A domain $\Omega\in\C^n$ is a Reinhardt domain  if for any $z=(z_1,\cdots, z_n)\in\Omega$, and any $(\theta_1,\cdots,\theta_n)\in\Rn$ we have $(e^{i\theta_1}z_1,\cdots,e^{i\theta_n}z_n)\in \Omega$. A Reinhardt domain $\Omega$ is called {\it complete} if for any $z=(z_1,\cdots,z_n)\in\Omega$ and for any $|\lambda_j|\leq 1, \lambda_j\in\C, j=1,\cdots,n$ one has $(\lambda_1z_1,\cdots,\lambda_nz_n)\in\Omega$.
\end{defi}

The pseudoconvexity and the hyperconvexity of a Reinhardt domain can be simply characterized using the logarithmic image $U:=\lbr (\log|z_1|,\cdots,\log|z_n|)\ |\ z\in\Omega\rbr$ and the axis set $M:=\{z\in\mathbb C^n| z_1\cdots z_n=0\}$.

\begin{prop}\label{pseudo}
 A Reinhardt domain in $\C^n$ is pseudoconvex if and only if $U$ is convex in $\Rn$ and moreover the following condition holds:

 if $M_j:=\lbr z\in M\ |\ z_j=0\rbr$ is such that $M_j\cap \Omega\neq \emptyset$ then for any $z\in\Omega$ the whole disc
 $$\lambda\longmapsto (z_1,\cdots, z_{j-1},\lambda z_j,z_{j+1},\cdots,z_n)$$
 ($|\lambda|\leq1$) belongs to $\Omega$.
\end{prop}
\begin{prop}\label{hyper}
 A pseudoconvex Reinhard domain $\Omega$ is hyperconvex if and only if it satisfies the the following property (known as the Fu condition):
 For every $j\in \{1,\cdots,n\}$ if $M_j\cap\overline{\Omega}\neq \emptyset$, then $M_j\cap\Omega\neq \emptyset$.
\end{prop}
More details on these domains can be found, for example, in the book \cite{JP08}.

For a fixed  Reinhardt domain $\Omega$ the rotationally invariant plurisubharmonic functions can be identified with convex functions on the logarithmic image through the logarithmic map. More precisely for a rotationally invariant (i.e. toric) plurisubharmonic function $\phi$, the function $u(x):=\phi(e^{x_1},\cdots,e^{x_n})$ is convex on $U$ - the logarithmic image of $\Omega$.


 From now on, we restrict our attention to the space $\tilde{\mathcal H}$ consisting of {\it toric} plurisubharmonic functions in $\mathcal H$. 
 \begin{lem}\label{convtor}Let $\varphi$ be a toric plurisubharmonic function from $\tilde{\mathcal H}$ and $u$ be as above. Then
	$$\phi(z_1,z_2,...,z_n)=u(\log|z_1|,\log|z_2|,,,\log|z_n|)=:u(x),$$
Moreover,
	\begin{equation*}
		\frac{\partial^2 \phi}{\partial z_i\partial \bar{z}_j}(z)=\frac{\partial}{\partial z_i}\left(\frac{\partial u}{\partial x_j}\frac{1}{2\overline{z}_j}\right)\\
		=\frac{\partial^2 u}{\partial x_i\partial x_j }(x)\frac{1}{ 4z_i\overline{z}_j}
	\end{equation*}

\end{lem}
From this Lemma, the following equality, holding whenever $u$ and $\phi$ are $C^2$ smooth, will be crucial:
\begin{equation}\label{torconv}
\frac{\partial^2 u}{\partial x_i\partial x_j }(x)=\frac{\partial^2 \phi}{\partial z_i\partial \bar{z}_j}(e^x)4e^{x_i}e^{x_j},
\end{equation}
where we used the suggestive notation $e^x=(e^{x_1},\cdots,e^{x_n})$.

Next proposition summarizes the findings in \cite{Ras17},\cite{AD21} and \cite{AD26} that will be relevant for us:

\begin{prop}\label{geodgeod} Let $\Omega$ be a bounded Reinhardt domain in $\mathbb C^n$ with $U$ - its  logarithmic image. Let also $\varphi$ be a weak geodesic joining $\varphi_0$ and $\varphi_1$ from $\tilde{\mathcal H}$. If $u_0(x):=\varphi_0(e^x),\ u_1(x):=\varphi_1(e^x)$, then
	$u(x,t):=\phi(e^x,e^t)$ is convex function on $U\times[0,1]$. Furthermore for any $(x,t)$ in $U\times[0,1]$ there are unique points $\xi=\xi(x,t)$ and $\eta=\eta(x,t)$ in $\overline{U}$ so that the line segment $l$ with endpoints $(\xi,0)$ and $(\eta,1)$ passes through $(x,t)$ and $u$ restricted to $l$ is affine.
\end{prop}

The affinity of $u$ along $l$ implies (see \cite{AD21}) that if both $\xi$ and $\eta$ are in $U$ then
\begin{equation}\label{dxideta}
 \forall j\in \{1,\cdots, n\}\quad D_{x_j}u_0(\xi)=D_{x_j}u_1(\eta),
\end{equation}
where $D_{x_j}$ denotes partial differentiation in the $j$-th direction.

As explained in \cite{AD26} the regularity of $\phi$ near $M\cap\Omega$ crucially depends on whether the vector $\xi(x,t)-\eta(x,t)$ is bounded. The following example, borrowed from \cite{AD26}, shows that this may fail if one of the functions is not strictly plurisubharmonic:

\begin{exa}\label{3rdexa}
Let $\Omega$ be the unit disc in $\mathbb C$ (hence the logarithmic image is the negative semi-axis). Let also $\varphi_0(z)=2(|z|^2-1),\ \varphi_1(z)=|z|^4-1$. Then the geodesic $\phi$ joining $\varphi_0$ and $ \varphi_1$ is equal to $(2-\log|\zeta|)[|z|^{\frac{4}{2-\log|\zeta|}}-1]$. The vectors $\xi$ and $\eta$ are given by the formula
$$\xi(x,t)=2\eta(x,t),\ \ \ \eta(x,t)=\frac{2x}{2-t}.$$

The difference $\xi-\eta$ is obviously unbounded as $x\rightarrow-\infty$.
\end{exa}

\section{Proof of the main theorem}

Note that existence of a function $\varphi\in \tilde{\mathcal H}$, that is, a toric plurisubharmonic function which is smooth strictly plurisubharmonic and vanishes on the boundary of $\Omega$ implies the hyperconvexity of $\Omega$. So by Proposition \ref{hyper}, $\Om$ satisfies the Fu condition. Hence, for a fixed $j\in \{1,\cdots,n\}$, we have the dichotomy:
\begin{itemize}
 \item {\bf either} $M_j\cap\overline{\Omega}=\emptyset$;
 \item {\bf or } $M_j\cap\Omega\neq\emptyset$.
\end{itemize}

In the first case, as $\overline{\Omega}$ is compact, the distance between $M_j$ and $\overline{\Omega}$ is positive. This implies that the logarithmic image $U$ of $\Omega$ is {\bf bounded} in the $j$-th coordinate direction i.e. the following holds: The set
$$U_j:=\lbr p\in\R\ |\ \exists\ {\bf p}=(p_1,\cdots,p_{j-1},p,p_{j+1},\cdots,p_n)\in U\rbr$$
is bounded. This is simply because $U_j$ is always bounded from above (as $\Omega$ is a bounded domain) whereas the positive distance yields the lower bound. As a result, for {\bf any} $\xi,\ \eta\in U$ the $j$-th coordinate of $\xi-\eta$ is uniformly bounded.

In the second case, as $\Omega$ is pseudoconvex, by Proposition \ref{pseudo}
we know that for any $z\in\Omega$ the whole disc
 $$\lambda\longmapsto (z_1,\cdots, z_{j-1},\lambda z_j,z_{j+1},\cdots,z_n)$$
 ($|\lambda|\leq1$) belongs to $\Omega$. This implies, in turn, on the logarithmic image $U$ the following holds:
  $x\in U$, and any $s<0$ the point $(x_1,\cdots, x_{j-1},x_j+s,x_{j+1},\cdots,x_n)$ belongs to $U$.

	Recall the equality $D_{x_j}u_0(\xi)=D_{x_j}u_1(\eta)$, which we rewrite as
	\begin{equation}\label{2dertrick}
		\int_{-\infty}^0D^2_{x_jx_j}u_0(\xi_1,\cdots, \xi_{j-1},\xi_j+s,\xi_{j+1},\cdots,\xi_n)ds
	\end{equation}
	$$=\int_{-\infty}^0D^2_{x_jx_j}u_1(\eta_1,\cdots, \eta_{j-1},\eta_j+s,\eta_{j+1},\cdots,\eta_n)ds,$$
where we crucially used the geometric property of the domain $U$ discussed above and coupled it with the convexity of $u_0$ and $u_1$ to justify the convergence of the integrals.

Now, due to Lemma \ref{convtor} and equation (\ref{torconv}) the equality can be further rewritten as
	\begin{equation}\label{2dertrickbis}
		\int_{-\infty}^0\frac{\partial^2\varphi_0}{\partial z_j\partial{\overline{z}_j}}(e^{\xi_1},\cdots, e^{\xi_{j-1}},e^{\xi_j+s},e^{\xi_{j+1}},\cdots,e^{\xi_n})4e^{2(\xi_j+s)}ds
	\end{equation}
	$$=\int_{-\infty}^0\frac{\partial^2\varphi_1}{\partial z_j\partial{\overline{z}_j}}(e^{\eta_1},\cdots, e^{\eta_{j-1}},e^{\eta_j+s},e^{\eta_{j+1}},\cdots,e^{\eta_n})4e^{2(\eta_j+s)}ds.$$

	Recall now that, by Step 1 of the proof of \cite[Theorem 5.3]{AD26}, both $e^\xi$ and $e^\eta$ are in a compact subset of $\Omega$ and hence the complex Hessians of $\varphi_0$ and $\varphi_1$ are uniformly bounded both from above and below along the integration paths. Hence for a constant $C$ dependent on $\varphi_0,\varphi_1$ one has
\begin{align*}
	&\int_{-\infty}^0\frac{\partial^2\varphi_0}{\partial z_j\partial{\overline{z}_j}}(e^{\xi_1},\cdots, e^{\xi_{j-1}},e^{\xi_j+s},e^{\xi_{j+1}},\cdots,e^{\xi_n})4e^{2(\xi_j+s)}ds\\
	&\qquad \geq C^{-1}\int_{-\infty}^0e^{2(\xi_j+s)}ds
= \frac{C^{-1}}2e^{2\xi_j},
\end{align*}
	while
	\begin{align*}
	&\int_{-\infty}^0\frac{\partial^2\varphi_1}{\partial z_j\partial{\overline{z}_j}}(e^{\eta_1},\cdots, e^{\eta_{j-1}},e^{\eta_j+s},e^{\eta_{j+1}},\cdots,e^{\eta_n})4e^{2(\eta_j+s)}ds\\
	&\qquad\leq C\int_{-\infty}^0e^{2(\eta_j+s)}ds
	= \frac{C}2e^{2\eta_j}.
	\end{align*}
Plugging these into (\ref{2dertrickbis}) one obtains
	$$e^{\xi_j}\leq C^2e^{\eta_j}.$$
Reversing the roles of $\varphi_0$ and $\varphi_1$ 	results in
	$$e^{\eta_j}\leq C^2e^{\xi_j}.$$
	This implies that the vector $\xi-\eta$ has $j$-th component  uniformly bounded as in the first case. This concludes the contention.

\end{document}